\documentclass
{amsart}

\usepackage[T2A]{fontenc}
\usepackage[cp1251]{inputenc}
\usepackage[english,russian]{babel}
\usepackage{amsmath}
\usepackage{amsfonts}
\usepackage{amssymb}

\usepackage{mathrsfs}

\usepackage{enumerate}
\usepackage{graphicx}

\usepackage[pdftex,unicode,colorlinks,linkcolor=blue,
citecolor=red,bookmarksopen,pdfhighlight=/N]{hyperref}

\numberwithin{equation}{subsection}

\usepackage[hyper]{amsbib}

\theoremstyle{plain}

\newtheorem{theoA}{Теорема A}
\newtheorem{maintheo}{Основная теорема}
\newtheorem{uniqtheo}{Теорема единственности}
\newtheorem{theoB}{Теорема B}
\newtheorem{theoS}{Теорема S}
\newtheorem{PLf}{Формула Пуанкаре\,--\,Лелона}
\newtheorem{lemma}{Лемма}[subsection]
\newtheorem{propos}{Предложение}[subsection]

\theoremstyle{definition}
\newtheorem{definition}{Определение}
\newtheorem{remark}{Замечание}[subsection]
\newtheorem{example}{Пример}[subsection]

\renewcommand{\leq}{\leqslant}
\renewcommand{\geq}{\geqslant}

\newcommand{\e}{\varepsilon}
\newcommand{\dd}{\,\mathrm{d}}

\newcommand{\const}{{\rm const}}
\newcommand{\rad}{{\rm rad}}

\def\RR{\mathbb R}
\def\CC{\mathbb C}
\def\NN{\mathbb N}
\def\BB{\mathbb B}
\def\SS{\mathbb S}

\DeclareMathOperator{\sbs}{sbs} 
\DeclareMathOperator{\clos}{clos}
\DeclareMathOperator{\Int}{int}
\DeclareMathOperator{\Har}{har}
\DeclareMathOperator{\Hol}{Hol}

\DeclareMathOperator{\Zero}{Zero}
\DeclareMathOperator{\sbh}{sbh}

\DeclareMathOperator{\supp}{supp}
\DeclareMathOperator{\Meas}{Meas}

\DeclareMathOperator{\dsbh}{\text{$\delta${\rm -sbh}}}
\DeclareMathOperator{\reg}{reg}

\begin{document}

\title{Нулевые множества голоморфных функций\\ в единичном шаре: \\ 
нерадиальные характеристики роста}
    
\author{Б.\,Н.~Хабибуллин, Ф.\,Б.~Хабибуллин}

\date{18.09.2018}

\maketitle


\begin{abstract}
Пусть $f$ --- ненулевая голоморфная функция в единичном шаре $\BB$ из $n$-ме\-р\-н\-о\-го комплексного евклидова пространства $\CC^n$, обращающаяся в нуль на множестве ${\sf Z}\subset \BB$ и  удовлетворяющая ограничению $|f|\leq \exp M$ на $\BB$, где $M\not\equiv \pm \infty$ --- $\delta$-субгармоническая в $\BB$ с зарядом Рисса $\mu_M$. Дается  шкала интегральных равномерных ограничений сверху на распределение множества ${\sf Z}$ через заряд $\nu_M$ в терминах $(2n-2)$-меры Хаусдорфа множества $\sf Z$, а также тестовых выпуклых радиальных функций и $\rho$-суб\-с\-ф\-е\-р\-и\-ч\-е\-с\-к\-их функций на единичной сфере $\SS\subset \CC^n$, которые при $n=1$ можно трактовать как $2\pi$-периодические  $\rho$-тригонометрически выпуклые функции на вещественной оси $\RR \subset \CC$.

Библиография: 25 названий.

Ключевые слова: голоморфная функция, нулевое множество, $\delta$-суб\-г\-а\-р\-м\-о\-н\-и\-ч\-е\-с\-к\-ая функция, заряд Рисса, $\rho$-тригонометрически выпуклая функция, $\rho$-субсферическая функция, мера Хаусдорфа   
\end{abstract}

\markright{Нулевые множества голоморфных функций в единичном шаре\dots}

\footnotetext[0]{Исследование выполнено за счет гранта Российского научного фонда (проект № 18-11-00002) --- первый автор, а также при финансовой поддержке РФФИ в рамках научного проекта  № 18-51-06002 --- второй автор}

\subsection{Введение.  О $\rho$-субсферических функциях}
Самые широкие применения  при исследовании поведения 
голоморфных и субгармонических функций на  {\it комплексной пло\-с\-к\-о\-с\-ти\/} $\CC$ и в угле  из $\CC$  находят {\it $\rho$-три\-г\-о\-н\-о\-м\-е\-т\-р\-и\-ч\-е\-с\-ки выпуклые функции\/} $h$ \cite{Levin56}--\cite[гл.~2]{GM} на связных подмножествах (интервалах) $I$ {\it вещественной оси\/} $\RR$  со значениями из  {\it расширенной вещественной оси\/}    
$\RR_{\pm \infty}:=\RR_{-\infty}\cup \RR_{+\infty}$, $\RR_{+\infty}:=\RR\cup \{+\infty \}$,
$\RR_{-\infty}:=\{-\infty\}\cup \RR$, $\RR^+:=\{x\in  \RR\colon x\geq 0\}$.
При $\rho \in \RR^+_*:=\RR^+\setminus \{0\}$ они полностью характеризуются неравенствами 
\begin{equation}\label{dftrc}
h(\theta)\leq \frac{\sin \rho (\theta_2-\theta)}{\sin \rho (\theta_2-\theta_1)}h(\theta_1)+\frac{\sin \rho (\theta-\theta_1)}{\sin \rho (\theta_2-\theta_1)}h(\theta_2)
\quad\text{\it для всех $\theta \in (\theta_1,\theta_2)\subset I$}
\end{equation}
{\it при\/} $0<\theta_2-\theta_1<\pi/\rho$. Функция   {\it $0$-три\-г\-о\-н\-о\-м\-е\-т\-р\-и\-ч\-е\-с\-ки выпуклая,\/} если она тождественная постоянная из $\RR_{\pm \infty}$. Далее удобно считать, что $\rho$-тригонометрически выпуклые функции не принимают значение $+\infty$, т.\,е. рассматриваются  $h\colon I \to \RR_{-\infty}$.

 Многомерные обобщения $2\pi$-периодических  $\rho$-три\-г\-о\-н\-о\-м\-е\-т\-р\-и\-ч\-е\-с\-ки выпуклых  фу\-н\-к\-ций $h\colon \RR \to \RR_{-\infty}$ --- это {\it $\rho$-субсферические функции\/} \cite[\S~4, определение 8]{Kon81},  \cite[\S~7, определение 10]{Kon84}, \cite[\S~1]{Kha91}, \cite[3, определение 3.1]{Kha99},   \cite[определение 4.2.1]{Khsur}, \cite[3.5, теорема S]{KhaAbdRoz18} на {\it единичной сфере\footnote{Метка-ссылка над знаками (не)равенства, включения или, более общ\'о, бинарного отношения и т.\,п. означает, что данное соотношение как-то связано с  отмеченной ссылкой.}\/} 
\begin{equation}\label{SB}
\SS:=\{ x\in \RR^m\colon |x|\overset{\eqref{norm}}{=}1\}, \quad \BB:=\{ x\in \RR^m\colon |x|\overset{\eqref{norm}}{<}1\},
\end{equation} 
из {\it $m$-мерного вещественного евклидова пространства\/} $\RR^m$,  $2\leq m\in \NN:=\{1,2,\dots\}$,  $\NN_0:=\{0\} \cup \NN$, со стандартной {\it евклидовой нормой\/}
\begin{equation}\label{norm}
 |x|:=\sqrt{x_1^2+x_2^2+\dots + x_m^2}, \quad x=(x_1,\dots x_m)\in \RR^m, 
\end{equation}  
или  {\it $n$-мерного комплексного евклидова пространства\/} $\CC^n$, $n\in \NN$, которое здесь часто можно и удобно  отождествлять с $\RR^{2n}$:
\begin{equation}\label{CR}
(x_1+iy_1,\dots , x_n+iy_n)\in \CC^n \; 
\longleftrightarrow\; (x_1,\dots , x_n,y_1,\dots,y_n)\in \RR^{2n}.
\end{equation}
Ранее классы $\rho$-субсферических функций, которые обозначаем далее через $\rho$-$\sbs (\SS)$,  применялись для распространения метода рядов Фурье для целых и мероморфных функций одной переменной из работ  Н.\,И.~Ахиезера,  Л.\,А.~Рубела и Б.\,А.~Тейлора \cite[\S~7]{GLO}  
 на субгармонические функции в $\RR^m$ в статьях   А.\,А.~Кондратюка 
\cite[\S~4]{Kon81}, \cite[\S~7]{Kon84}. 
Кроме того, они были использованы для исследования   нулевых множеств целых функций с ограничениями на их рост  в $\CC^n$ и двойственных им условий полноты экспоненциальных систем в пространствах голоморфных функций в областях из $\CC^n$ в работах первого из авторов \cite[\S~4]{Kha91}, \cite[теоремы 2.1, 3.1, 3.2]{Kha99}, \cite[теоремы 3.3.5, 4.2.7]{Khsur}, 
 \cite[5.1.6]{KhaAbdRoz18}.  Прямые применения классов $\rho$-$\sbs (\SS)$ к исследованию поведения голоморфных функций в  областях из $\CC^n$, отличных от $\CC^n$, при  $n>1$ нам неизвестны. Здесь мы используем  классы $\rho$-$\sbs (\SS)$ для исследования распределения нулевых множеств голоморфных функций с ограничениями на их рост в {\it единичном шаре\/}  
\begin{equation}\label{DSd}
\BB\overset{\eqref{SB}}{:=}\bigl\{z\in \CC^n \colon |z| < 1 \bigr\}\subset \CC^n.
\end{equation}
Если для $\SS$  и $\BB$ в  $\RR^m$ или  в $\CC^n$, отождествленном с $\RR^{2n-1}$ как в  \eqref{CR}, необходимо указать  размерность, то пишем соответственно $\SS_{m-1}$ и 
 $\BB_m$ или $\SS_{2n-1}$ и $\BB_{2n}$, где нижний индекс указывает на <<вещественную>> размерность.  В обозначениях  $\Gamma$ для гамма-функции и $a^+:=\max \{0,a\}$ для положительной части $a$ полагаем 
\begin{subequations}\label{bsBS}
\begin{align}
b_m&:=\frac{\pi^{m/2}}{\Gamma (m/2+1)}=
\begin{cases}
{\pi^n}/{n!} &m=2n\in 2\NN,\\
\dfrac{n! \, 2^{2n+1}\pi^n}{(2n+1)!}&m=2n+1\in 2\NN+1,
\end{cases} , \quad s_{m-1}:=mb_m 
\tag{\ref{bsBS}b}\label{df:sbm} 
\\
\intertext{ --- соответственно объем шара $\BB_m\subset \RR^m$ и площадь сферы $\SS_{m-1}\subset \RR^m$,}
d_{m-1}&:=\bigl(1+(m-3)^+\bigr)\,s_{m-1} , \quad 
d_{2n-1}\overset{\eqref{df:sbm}}{=} \frac{2\pi^n \max\{1,2n-2\}}{(n-1)!}
\tag{\ref{bsBS}d}\label{df:sbm+}
\end{align}
\end{subequations}
--- некоторые необходимые далее нормирующие\footnote{В \cite[(1.2), (1.4)]{KhaRoz18}
числа $d_{m-1}$ и $d_{2n-1}$ некорректно трактуются как площади соответственно  $s_{m-1}$ и $s_{2n-1}$ сферы $\SS$, а в формуле для $b_{2n-2}=\pi^{n-1}/(n-1)!$ справа присутствует лишний множитель $\max\{1,2(n-2) \}$. Впрочем, нормировки мер Рисса и Хаусдорфа в \cite[(1.2), (1.3)]{KhaRoz18} описаны верно, поэтому  указанные неточности и описки не повлияли на справедливость 
результатов из \cite{KhaRoz18}.} множители\,/\,делители  \cite[гл. 1, \S~2]{Ron71}.

Через  $\Hol(S)$ при $S\subset \CC^n$ и  $\Har(S)$, $\sbh (S)$, $\dsbh(S):=\sbh (S)-\sbh (S)$, $C^k(S)$ для  $k\in \NN\cup \{\infty\}$ при $S\subset \RR^m$ 
обозначаем соответственно классы  {\it голоморфных,  гармонических, субгармонических \cite{HK}, $\delta$-субгармонических\/ {\cite[3.1]{KhaRoz18}--\cite{Gr}},  непрерывно $k$ раз дифференцируемых\/} функций $v$ на некотором открытом множестве $\mathcal O_v\supset S$ --- своем для каждой функции $v$. Но $C(S)$ --- класс всех {\it непрерывных\/} функций на $S$ в топологии, индуцированной с $\CC^n$ или $\RR^m$. Функции, тождественно равные $-\infty$ или $+\infty$ на  $S$ обозначаем соответственно как $\boldsymbol{-\infty}$ или  $\boldsymbol{+\infty}$. Кроме того, 
полагаем 
\begin{equation*}
\sbh_*(S):=\sbh(S)\setminus \{\boldsymbol{-\infty} \}, \; 
\dsbh_*(S):=\dsbh(S)\setminus \{\boldsymbol{\pm\infty} \},  
\; \Hol_*(S):=\Hol(S)\setminus \{0\}.
\end{equation*}
При этом  одним и тем же символом $0$ обозначаем, по контексту, число нуль, начало отсчета, нулевой вектор, нулевую функцию, нулевую меру и т.\,п.  {\it Положительность\/} всюду понимается как $\geq 0$; противоположное $\leq 0$ --- отрицательность.  

\begin{definition}[{\rm \cite[определение 1.2]{Kha91}, \cite[определение 2]{KhaAbdRoz18}}]\label{df:Krr}
Пусть  $x\cdot y$ --- скалярное произведение радиус-векторов в $\RR^m$ точек $x,y\in \SS_{m-1}$,
 $\rho\in \RR^+_*$ и $0<r\leq 1$. {\it Ядром усреднения\/} $K_{\rho,r}\colon \SS_{m-1}\times \SS_{m-1}\to \RR^+$ для субсферичности называем положительную функцию, которая при  
$x\cdot y>\sqrt{1-r^2}$ определена как
\begin{multline}\label{Krr}
K_{\rho,r}(x, y):=\frac{1}{\rho+m}\left(\Bigl(x\cdot y +\sqrt{(x\cdot y)^2-(1-r^2)}\Bigr)^{\rho+m}\right.\\
\left.-\Bigl(x\cdot y -\sqrt{(x\cdot y)^2-(1-r^2)}\Bigr)^{\rho+m}\right)\\
=\frac{1}{\rho+m}\left(\Bigl(\cos \varphi  +\sqrt{r^2-\sin^2 \varphi}\,\Bigr)^{\rho+m}
-\Bigl(\cos \varphi  -\sqrt{r^2-\sin^2 \varphi}\,\Bigr)^{\rho+m}\right),
\end{multline}
где $\varphi:=\arccos (x\cdot y)$ --- это угол между радиус-векторами точек $x,y\in \SS_{m-1}$, а при 
$x\cdot y\leq \sqrt{1-r^2}$ полагаем $K_{\rho,r}(x, y):=0$.
\end{definition}

При $\rho=0$ по определению класс $0$-$\sbs(\SS)$ --- это тождественные постоянные
из $\RR_{-\infty}$. Определяет  $\rho$-субсферические функции $h\neq \boldsymbol{-\infty} $ при $\rho>0$ 
\begin{theoS}[{\rm \cite[предложение 1.3]{Kha91}, 
\cite[теорема S]{KhaAbdRoz18}}] Пусть $h\colon \SS_{m-1} \to \RR_{-\infty}$ --- функция на единичной сфере $\SS_{m-1}\subset \RR^m$ и $h\neq \boldsymbol{-\infty}$, $\rho\in \RR^+_*$.  Следующие четыре утверждения эквивалентны:
\begin{enumerate}[{\rm (s1)}]
\item\label{ss1} $h$ --- субсферическая функция порядка $\rho$, т.\,е. из  класса $\rho$-$\sbs (\SS_{m-1})$;
\item\label{ss4} функция $h$ полунепрерывна сверху, интегрируема по мере $\sigma_{m-1}$ на $\SS_{m-1}$ и
\begin{equation}\label{Lrho}
\mathcal L_{\rho} h\geq 0, \quad   \mathcal L_{\rho} :=\Delta_{\SS} +\rho (\rho+m-2),
\end{equation}
где $\Delta_{\SS}$ --- сферическая часть на $\SS_{m-1}$ оператора Лапласа $\Delta$ на $\RR^m$
\cite[гл.~1, \S~7]{Ron71}, или оператор Лапласа\,--\,Бельтрами на $\SS_{m-1}$, а  неравенство из \eqref{Lrho}
выполнено в пространстве $\mathfrak D'(\SS_{m-1})$ распределений, или обобщенных функций, на $\SS_{m-1}$, и из этого неравенства следует, что $\mathcal L_{\rho} h$ совпадает в  $\mathfrak D'(\SS_{m-1})$ с некоторой   борелевской регулярной  положительной мерой на $\SS_{m-1}$
\cite[определение 1.1]{Kha91}, \cite[определение 3.1]{Kha99}, \cite[определение 4.2.1]{Khsur};
\item\label{ss2} функция $H(x)=h\bigl(x/|x|\bigr)|x|^{\rho}$, $x\in \RR^m\setminus \{0\}$,  доопределенная нулем в точке $x=0$, т.\,е. при $H(0):=0$,  субгармоническая  в $\RR^m$;
\item\label{ss3} функция $h$ полунепрерывна сверху и для любой точки $x\in \SS_{m-1}$ найдётся число $r_x\in (0,1)$, для которого 
\begin{equation}\label{hssf}
h(x)\leq \frac{1}{b_m r^m} \int_{\SS_{m-1}} h(y)K_{\rho,r}(x, y) \dd \sigma_{m-1}(y)
\quad \text{при всех\/  $0<r\leq r_x$},
\end{equation}
где $\dd \sigma_{m-1}$ --- элемент площади на $\SS_{m-1}$, т.\,е. $\sigma_{m-1}$  --- это $(m-1)$-мера Хаусдорфа на $\RR^m$ по  определению\/ \eqref{df:sbm}.
\end{enumerate} 
\end{theoS}
При $m=2$, или  для $\CC$, отождествленного с  $\RR^2$ как в \eqref{CR}, 
функция $h$  принадлежит классу $\rho$-$\sbs (\SS_1)$ тогда и только тогда,  когда функция
 $h(e^{i\theta})$, $\theta \in \RR$, --- $2\pi$-периодическая $\rho$-тригонометрически выпуклая функция на $\RR$. Следует заметить, что соотношение \eqref{hssf} не дает прямо неравенство    \eqref{dftrc}, хотя утверждение  (s\ref{ss3}) теоремы S эквивалентно свойству \eqref{dftrc} \cite[теорема 60]{GM} \cite[определение 1.3]{Kha91}, \cite[определение 3.1]{Kha99}. 

Класс функций $\rho$-$\sbs (\SS)$ замкнут относительно операции поточечного максимума, что легко следует из теоремы S(s\ref{ss2}). В частности, 
\begin{enumerate}[{\rm (i)}]
\item\label{tr:ii}  если  $h\in \rho\text{-}{\sbs} (\SS)$, то функция 
$h^+:=\max\{0,h\}$ принадлежит классу 
\begin{equation}\label{SS+}
 \rho\text{-}\sbs^+ (\SS):=
\{ h\in \rho\text{-}\sbs (\SS)\colon h\geq 0 \text{ на $\SS$}\}. 
\end{equation}
\item\label{tr:vi}  Если  $0\leq \rho\leq \rho'\in \RR^+$, то по теореме S(s\ref{ss4}) имеем  
$\rho\text{-}\sbs^+ (\SS) \overset{\eqref{SS+}}{\subset} \rho'\text{-}\sbs^+(\SS)$.  
\item\label{tr:vii} Если последовательность функций 
$h_n\in \rho\text{-}\sbs (\SS)$, $n\in \NN$, убывает, то по теореме S(s\ref{ss3}) поточечный предел
\begin{equation}\label{limh}
h:=\lim_{n\to \infty}h_n
\end{equation}
 принадлежит тому же классу $\rho$-$\sbs (\SS)$. Здесь класс $\rho$-$\sbs (\SS)$ можно заменить на класс $\rho$-$\sbs^+ (\SS)$ из \eqref{SS+}.
\item\label{tr:6} Если $h\in \rho\text{-}\sbs (\SS)$, то для любого $k\in \NN$ найдется убывающая последовательность функций  $h_n\in \rho\text{-}\sbs (\SS)\cap C^k(\SS)$, $n\in \NN$, для которой выполнено равенство \eqref{limh} \cite[предложение 1.7]{Kha91}.
\end{enumerate}

\subsection{Меры Хаусдорфа и дивизоры нулей} 
Для $p\in \RR^+$ через $\sigma_p$ обозначаем {\it $p$-мерную (внешнюю) меру Хаусдорфа,\/} или {\it $p$-меру Хаусдорфа,\/} в  $\RR^m$.  
 В настоящей статье  $p$-мера Хаусдорфа используется лишь при  целом $p{\in} \NN_0$:
\begin{equation}\label{df:spb}
\sigma_p(S)\overset{\eqref{df:sbm}}{:=}b_p \lim_{0<r\to 0} \inf \biggl\{\sum_{j\in \NN}r_j^p\,\colon\,,  
S\subset \bigcup_{j\in \NN}B(x_j,r_j), \; 0\leq r_j<r\biggr\}, 
\end{equation}
где $B(x,r):=x+r\BB$ --- {\it  открытый шар в\/  $\RR^m$ с центром $x\in \RR^m$ радиуса\/ $r$.}  В такой нормировке при $p=0$ для любого подмножества $S\subset  \RR^m$   его $0$-мера Хаусдорфа $\sigma_0(S)$ равна мощности $S$, т.\,е. числу точек в $S$, а при $p=m$ мера   $\sigma_m$ --- мера Лебега в $\RR^m$. 
В $\CC^n$, отождествленном с $\RR^{2n}$ как в \eqref{CR}, в настоящей работе будем использовать в подавляющем числе случаев лишь  $(2n-2)$-меру Хаусдорфа $\sigma_{2n-2}$.

Пусть $D$ --- область в $\CC^n$. Следуя \rm \cite[гл.~4]{Khsur}, \cite[\S~11]{Ron77}--\cite[гл.~1]{Chi},  {\it дивизором нулей\/} функции $f\in \Hol_*(D)$ называем функцию  $\Zero_f\colon D \to \NN_0$, равную кратности нуля функции $f$ в каждой точке $z\in D$. 
  Для $f=0\in \Hol (D)$ по определению $\Zero_0= \boldsymbol{+\infty}$ на  $D$. Дивизор нулей $\Zero_f$ полунепрерывен сверху в $D$.  Носитель $\supp \Zero_f$ --- главное аналитическое множество  чистой размерности $n-1$ над $\CC$  и размерности $2n-2$ над $\RR$, для которого $\reg \supp \Zero_f$ --- множество регулярных точек, и всегда  $\sigma_{2n-2} (\supp \Zero_f\setminus \reg \supp \Zero_f)=0$. Пусть   $\reg \supp \Zero_f=\cup_j {\sf Z}_j$ --- представление в виде объединения не более чем счётного числа связных компонент ${\sf Z}_j$, $j=1,2, \dots$. Тогда семейство $\{{\sf Z}_j\}$   локально конечно в $D$, т.\,е. каждое подмножество $S\Subset D$ пересекается лишь с конечным числом компонент ${\sf Z}_j$. Дивизор нулей  $\Zero_f$ постоянен на каждой компоненте ${\sf Z}_j$, т.\,е. однозначно определено значение $\Zero_f({\sf Z}_j)$ для каждого $j\in \NN$.    
  Каждому дивизору нулей $\Zero_f$  сопоставляем  положительную {\it считающую меру нулей\/}   $n_{\Zero_f}$, 
определяемую  как мера Радона равенствами 
\begin{equation}\label{nZ}
n_{\Zero_f}(\varphi ):=:\int \varphi  \dd  n_{\Zero_f} \overset{\eqref{df:spb}}{:=}\int \varphi \Zero_f \dd  \sigma_{2n-2} 
\end{equation}по всем финитным функциям $\varphi\in C(D)$, или эквивалентно,  как  борелевская мера на $D$ по правилу 
\begin{equation}\label{nZB}
n_{\Zero_f}(B)\overset{\eqref{df:spb}}{=}\sum\limits_{j} \Zero_f({\sf Z}_j)  \sigma_{2n-2}(B\cap {\sf Z}_j)
\end{equation} 
для любых борелевских подмножеств $B\subset  D$.

\begin{PLf}[{\cite[1.2.4]{KhaRoz18}}]
Пусть $n\in \NN$, $D\neq \varnothing$ --- область в $\CC^n$,  $f\in \Hol_* (D)$. Для   меры Рисса $\mu_{\log |f|}$ функции $\log |f|\in \sbh_*(D)$ имеем равенства
\begin{equation}\label{nufZ}
\mu_{\log |f|}{:=} \frac{1}{d_{2n-1}}\Delta \log |f| \overset{\eqref{df:sbm}}{=}
\frac{(n-1)!}{2\pi^n \max\{1,2n-2\}}\Delta \log |f| {=}n_{\Zero_f}.
\end{equation}
\end{PLf}
Функция  ${\rm Z}\colon D\to \RR^+$  ---  {\it поддивизор нулей  для\/} $f\in \Hol (D)$, или {\it поддивизор дивизора нулей\/}  $\Zero_f$,  если  ${\rm Z}\leq \Zero_f$ на $D$. Очевидно,  для $f\in  \Hol (D)$ ее дивизор нулей $\Zero_f$  --- это и ее  поддивизор нулей.
  
\subsection{Основные результаты}\label{sec1}

Далее $\Meas\, (S)$  --- класс {\it борелевских вещественных мер,\/} или {\it зарядов,\/}  на   $S\subset \CC_{\infty}$;  ${\Meas}^+ (S)\subset {\Meas\,} (S)$ --- подкласс положительных мер в ${\Meas\,} (S)$.  Интегралы по положительным мерам $\mu \in \Meas^+(S)$
всюду, вообще говоря, понимаем как {\it верхние интегралы\/} \cite{B} с естественным продолжением на интегралы по зарядам $\mu \in \Meas\,(S)$. Для заряда  $\mu\in \Meas\, (S)$ его сужение на множество $X$ обозначаем как $\mu\bigm|_X$.
 Для $S\subset \RR^m$ {\it меру (заряд) Рисса\/} функции $u\in \sbh_*(S)$ (соответственно $u\in \dsbh_*(S)$)  обозначаем  как $\mu_{u}\overset{\eqref{df:sbm+}}{:=}\frac{1}{d_{m-1}}\Delta u\in {\Meas} ^+(S)$ (соответственно $\in {\Meas} (S)$). 
 
\begin{definition}[{\cite[(3.1)]{Kha91}, \cite[(0.2)]{Kha99}}]\label{df:mht}
{\it Радиальную  считающую функцию   заряда\/ $\mu \in \Meas\, (\BB)$ с весом\/}  $h\colon \SS\to \RR_{\pm \infty}$ определяем как функцию 
\begin{equation}\label{ntk}
\mu^{\rad}(r;h):=\int_{r\BB} h \bigl(x/|x|\bigr) \dd \mu (x), \quad r\in (0,1),
\quad  r\BB:=\{rx\colon x\in \BB \}.
\end{equation} 
При $h\equiv 1$ полагаем $\mu^{\rad}(r):=\mu^{\rad}(r;1)$ --- классическая {\it радиальная считающая функция заряда\/ $\mu \in \Meas\,(\BB)$.} 
{\it Радиальную считающую функцию  функции\/  
$\sf Z \colon \BB \to \RR$ с весом\/}  $h\colon \SS\to \RR_{\pm\infty}$  {\it относительно $p$-меры Хаусдорфа\/} определяем как функцию 
\begin{equation}\label{rcfZ}
{\sf Z}_p^{\rad}(r ;h):=\int_{r\BB} {\sf Z} (z)  h \bigl(z/|z|\bigr) \dd \sigma_p(z), \quad  r\in (0,1).
\end{equation}
При  $\BB\subset \CC^n$ и   $p=2n-2$ индекс $p$ в ${\sf Z}_p^{\rad}(\cdot ;h)$ из \eqref{rcfZ} опускаем и пишем ${\sf Z}^{\rad}(\cdot ;h)$. 
\end{definition}

Для  $-\infty\leq r<R\leq +\infty$  далее всюду  
\begin{equation}\label{irR}
\int_r^R\ldots :=\int_{(r,R)} \ldots \; .
\end{equation} 

\begin{uniqtheo}[{\rm (индивидуальная, $\BB\overset{\eqref{DSd}}{\subset} \CC^n$})]\label{th:ind}
Пусть   $M\in \dsbh_*(\BB) $ --- функция  с зарядом Рисса  $\mu_M \overset{\eqref{df:sbm}}{:=}\frac{1}{d_{2n-1}}\Delta M$, $f\in \Hol(\BB)$  и  $|f|\leq \exp M$ на $\mathbb B$. Допустим, что ${\sf Z}$ --- поддивизор дивизора нулей $\Zero_f$. 
Пусть  $g\colon \RR^+ \to \RR^+$ --- выпуклая функция с  $g(0)=0$ и 
$h\in \rho\text{-}\sbh^+ (\SS)$ для некоторого $\rho \in \RR^+$. Если  
\begin{equation}\label{1.n}
\int_{1/2}^{1} g\bigl(2^{2n-1}(2n-1)(1-r)\bigr) {\dd}\mu_M(r;h)<+\infty,
\end{equation} 
но 
\begin{equation}\label{cuZ}
\int_{1/2}^{1} g(1-r) {\dd} {\sf Z}^{\rad}(r ;h)=+\infty, 
\end{equation}
то  $f=0$.
\end{uniqtheo}
\begin{remark}\label{rem1}
В случае  $M=0$ с  $\mu_M=0$, $g(x)\equiv x^p$, $x\in \RR^+$,   $h\equiv 1\in 0\text{-}\sbs^+(\SS)$,  условие  \eqref{cuZ} противоречит условию Бляшке
 \begin{equation}\label{NHD}
\int_{{\sf Z} \cap (\BB\setminus \frac12 \BB)} (1-|z|) {\dd} \sigma_{2n-2}(z) <+\infty, \quad \e\in (0,1). 
\end{equation}
Таким образом, многомерная версия теоремы Неванлинны --- критерий 
 Г.~М.~Хенкина\,--\,Х.~Скоды о распределении  нулевого множества
ограниченной голоморфной функции   \cite[6.5]{Hen85}  --- показывает, что наша теорема единственности в этом случае оптимальна. Аналогично, для субгармонической в $\BB$ функции 
\begin{equation}\label{Mp}
M (z)= \const\frac{1}{(1-|z|)^p}\, ,   \quad z\in \BB, \quad p\in \RR_*^+,
\end{equation}
из $|f|\leq \exp M$ на $\BB$ и  ${\sf Z}\leq \Zero_f$  по теореме-критерию  Ш.~М.~Даутова\,--\,Х.~Скоды о нулевых множествах для весовых пространств Джрбашяна\,--\,Неванлинны \cite[6.5]{Hen85}   следует соотношение 
\begin{equation}\label{NHDp}
\int_{{\sf Z} \cap (\BB\setminus \frac12 \BB)} (1-|z|)^{p+1+\e} {\dd} \sigma_{2n-2}(z) <+\infty \quad \text{для любого $\e\in \RR_*^+$,} 
\end{equation}
что вполне согласуется с нашей теоремой  единственности.
\end{remark}
Через  $\const_{a_1,a_2,\dots} \in \RR$ обозначаем постоянные, которые зависят от  $a_1,a_2,\dots$ и, если не оговорено противное, только от них; $\const^+_{\dots} \geq 0$. 

\begin{maintheo}
[{\rm (равномерная,  $\BB\overset{\eqref{SB}}{\subset} \RR^m$)}] 
Пусть две функции $u\in \sbh_*(\BB)$ и  $M\in \dsbh_*(\BB)$ соответственно с мерой Рисса $\mu_u:=\frac{1}{d_{m-1}}\Delta u\in \Meas^+(\BB)$ и зарядом Рисса   $\mu_M:=\frac{1}{d_{m-1}}\Delta M\in \Meas (\BB)$ удовлетворяют неравенству $u\leq M$ на $\BB$. Пусть  $\rho\in \RR^+$. Тогда существует постоянная  $C:=\const^+_{\rho, M, u}\geq 0$  для которой неравенство 
\begin{equation}\label{uM}
\int_{1/2}^{1} g\Bigl(\frac{1}{r^{m-1}}-1\Bigr) \dd \mu^{\rad}_u(r;h)\overset{\eqref{ntk}}{\leq} 
\int_{1/2}^{1} g\Bigl(\frac{1}{r^{m-1}}-1\Bigr) \dd \mu^{\rad}_M(r;h)+C
\end{equation}
выполнено при  любых 
\begin{enumerate}
\item[{\rm [g]}]  выпуклой функции  $g\colon \RR^+\to \RR^+$ с  $g(0)=0$ и $g(2^{m-1}-1)\leq 1$,
\item[{\rm [h]}] $\rho$-субсферической функции $h\colon \SS\to [0,1]$.
\end{enumerate}
В частности, если\/ $\sf Z$ --- поддивизор нулей для какой-нибудь функции\/ 
$f\in \Hol_*(\BB)$, $\BB\overset{\eqref{DSd}}{\subset}\CC^n$, удовлетворяющей неравенству\/ $|f|\leq \exp M$ на $\BB$,  то для некоторой постоянной $C:=\const^+_{n,\rho,M,{\sf Z}}$ имет место неравенства 
\begin{equation}\label{uMf}
\int_{1/2}^{1} g\Bigl(\frac{1}{r^{2n-1}}-1\Bigr) {\dd} {\sf Z}^{\rad}(r ;h)
 \overset{\eqref{rcfZ}}{\leq} 
\int_{1/2}^{1} g\Bigl(\frac{1}{r^{2n-1}}-1\Bigr) \dd \mu^{\rad}_M(r;h)+C
\end{equation}
при любых функциях $g$ и $h$ из\/ {\rm [g]--[h]} с $m=2n$.
\end{maintheo}
\begin{remark}\label{rem2}
При выборе $M:=\log |f|$ в заключительной части основной теоремы в \eqref{uMf} получаем равенство с $C=0$. Таким образом, основная теорема точна с точностью до аддитивного слагаемого $C$. 
\end{remark}
\subsection{Примеры $\rho$-субсферических функций}
Пусть  $\rho\in \RR^+$, $\SS=\SS_{m-1}\subset \RR^m$.
\begin{example}\label{ex1}   Пусть  $s_0\in \SS$, $\angle (s,s_0)\in [-\pi,\pi]$ --- угол между радиус-векторами  точек  $s_0$  и $s\in \SS$.  Функция 
\begin{equation*}
h(s):=\begin{cases}
\cos \rho \angle (s,s_0)\quad &\text{при  $\bigl|\angle (s,s_0)\bigr|<\frac{\pi}{2\rho}$},\\
0\quad &\text{при  $\bigl|\angle (s,s_0)\bigr|\geq\frac{\pi}{2\rho}$},
\end{cases}
\qquad s\in \SS,
\end{equation*}  
принадлежит классу  $\rho$-$\sbs^+(\SS)$.
\end{example}
\begin{example}\label{ex2}
Ядро усреднения $K_{\rho,r}\colon \SS_{m-1}\times \SS_{m-1}\to \RR^+$ из определения \ref{df:Krr}, заданное  равенством \eqref{Krr},  при фиксации одного из двух аргументов  
становится положительной $\rho$-субсферической функцией по другой переменной. 
\end{example}

\begin{example}\label{ex3}  Пусть  $S\subset \RR^m$ --- ограниченное множество. 
Сужение на единичную сферу $\SS$ опорной функции множества $S$, определенное  как
\cite[определение 3.8]{Maergoiz}, \cite{GL}
\begin{equation*}
K_{S}(s):=\sup_{s'\in S}  (s\cdot s') , \quad s\in \SS,
\end{equation*}
принадлежит классу  $1$-$\sbs(\SS)$. Если $0\in S$, то  $K_{S}\in 1\text{-}\sbs^+(\SS)$.  
\end{example}

\begin{example}\label{ex4}  Пусть $u\in \sbh_*(\RR^m)$ или $u=\log |f|$, где $f\in \Hol_*(\CC^n)$. Если
\begin{equation*}
\limsup_{x\to \infty}\frac{u(x)}{|x|^{\rho}}<+\infty,
\end{equation*}
то сужение на $\SS$ полунепрерывной сверху регуляризации  радиального  $\rho$-индикатора 
\begin{equation*}
H_u(x):=\limsup_{r\to +\infty}\frac{u(rx)}{r^{\rho}}\, , \quad x\in \RR^m, \quad H(0):=0,
\end{equation*}
принадлежит классу $\rho$-$\sbs (\SS)$   \cite[гл.~3, \S~5]{Ron71}, \cite[1.3]{Ron86}, \cite{GL},  \cite{Ron92},  а те же операции применительно к $H_u^+$ дают функцию из $\rho$-$\sbs^+ (\SS)$. 
\end{example}
\begin{remark}
На основе базовых примеров \ref{ex1}--\ref{ex4} с помощью операций сложения, умножения на положительное число, взятия точной верхней грани ограниченного сверху семейства $\rho$-субсферических функций с последующей полунепрерывной сверху регуляризации, а также многих других действий, основанных на теореме S и сохраняющих $\rho$-субсферичность, можно строит разнообразные виды функций класса $\rho$-$\sbs(\SS)$  и $\rho$-$\sbs^+(\SS)$.  
\end{remark}

\subsection{Субгармонические тестовые функции и их роль}\label{ssec1_2} 

$\CC_{\infty}^n:=\CC^n\cup  {\infty}$  и $\RR^m_{\infty}:=\RR^m\cup \{\infty\}$ --- одноточечные компактификации Александрова соответственно $\CC^n$ и $\RR^m$.  
 Далее те определения или понятия, которые в рамках соответствия \eqref{CR} сразу переносятся с $\RR^{2n}_{\infty}$  на $\CC^n_{\infty}$, описываем и определяем только для $\RR^m_{\infty}$.
Для подмножества  $S\subset \RR^m_{\infty}$ 
через  $\clos S$, $\Int S$ и  $\partial S$ обозначаем 
{\it замыкание, внутренность и границу\/} $S$  в $\RR^m_{\infty }$. Открытое связное подмножество в  $\RR^m_{\infty}$ --- {\it  (под)область\/} в  $\RR^m_{\infty}$.  Для $S_0\subset S\subset  \RR^m_{\infty}$  пишем $S_0\Subset S$, если $\clos\, S_0$ --- компакт в  $S$ с топологией, индуцированной  на  $S$  с $\mathbb \RR^{m}_{\infty}$.  Пусть  
\begin{equation}\label{SD0D}
\varnothing \neq S\Subset D\subset \RR^m_{\infty},\quad\text{где  $D\neq \RR^m_{\infty}$ {\it --- область}.} 
\end{equation} 
Для функции $v\colon D\setminus S\to \RR$ пишем 
\begin{equation}\label{dD}
\lim_{\partial D}v=0,\quad\text{\it если\/  $\lim_{D\ni z'\to z}v(z')=0$ для всех\/ $z\in \partial D$.}  
\end{equation}
По определению положим 
\begin{subequations}\label{sbh}
\begin{align}
\sbh_0(D\setminus S) &:=\Bigl\{v\in \sbh(D\setminus S)\colon 
 \lim_{\partial D}v\overset{\eqref{dD}}{:=}0\Bigr\},
\tag{\ref{sbh}o}\label{sbh0}\\
\sbh_0^+(D\setminus S)&\overset{\eqref{sbh0}}{:=}\bigl\{v\in \sbh_0(D\setminus S)\colon v\geq 0\text{ on }D \bigr\}. 
\tag{\ref{sbh}+}\label{sbh0+}
\end{align}
\end{subequations}

\begin{definition}[{\rm \cite[определение 1]{KhaRoz18}}]\label{def:tesf} 
Функцию  $v\overset{\eqref{sbh0+}}{\in} \sbh^+_0 (D\setminus S)$ называем  положительной {\it субгармонической тестовой функцией на\/  $D$ вне\/} $S$, если функция  $v$ ограничена в $D\setminus S$. Класс таких функций $v$ будем обозначать как $\sbh_0^{+}(D\setminus S; <+\infty)$. Для  $b\in \RR^+$ полагаем 
\begin{equation}\label{<b}
\sbh_0^{+}(D\setminus S; \leq b)\overset{\eqref{sbh0+}}{:=}\Bigl\{v\in \sbh_0^{+}(D\setminus S;<+\infty)\colon \sup_{D\setminus S} v\leq b\Bigr\}.
\end{equation}
Таким образом,
\begin{equation*}
\sbh_0^{+}(D\setminus S; <+\infty)=\bigcup_{b\in \RR^+}\sbh_0^{+}(D\setminus S; \leq b).
\end{equation*}
\end{definition}
Основную роль будет играть следующая
\begin{theoA}[{\rm (\cite[основная теорема]{KhaRoz18}, см. и \cite{KhaAbdRoz18},
\cite{KhaTam17_A}, \cite{KhaTam17_L})}]  
Пусть $ M\in \dsbh_* (D)$ --- функция с зарядом Рисса $\mu_M=\frac{1}{d_{m-1}}\Delta M\in \Meas (D)$ и   
\begin{equation}\label{SD0D+}
\varnothing \neq \Int S\subset S=\clos S \overset{\eqref{SD0D}}{\Subset} D\subset \CC_{\infty}\neq  D.
\end{equation}
Тогда для любой точки  $x_0\in \Int S$ с $M(x_0)\in  \RR$, любого числа  $b{\in} \RR_*^+$,  любой регулярной для задачи Дирихле\/ {\rm \cite[2.6]{HK}}  области  $\widetilde{D}\subset \CC_{\infty}$ с функцией Грина $g_{\widetilde{D}}(\cdot , x_0)$ с полюсом в $x_0$ при предположении    $S\Subset \widetilde{D}\subset D$  и  $\CC_{\infty}\setminus \clos \widetilde{D}\neq \varnothing$, любой функции $u\in \sbh_* (D)$ с мерой Рисса $\mu_u=\frac{1}{d_{m-1}}\Delta u\in \Meas^+(D)$ и ограничением $u\leq M$ на  $D$, а также любой субгармонической тестовой функции   $v\overset{\eqref{<b}}{\in}  \sbh_0^+(D\setminus S;\leq b) $ выполнено неравенство 
\begin{equation}\label{mest+}
\widetilde{C} u(x_0) 	+\int_{D\setminus S}  v \dd {\mu}_u	\leq \int_{D\setminus S}  v \dd {\mu}_M	+\int_{\widetilde{D}\setminus S} v \dd {\mu}_M^-   +\widetilde{C}\, \overline{C}_M,
\end{equation}
где  $\mu_M^-$ --- нижняя вариация заряда $\mu_M$,  $\widetilde C, \overline C_M$ --- постоянные, определенные как  
\begin{subequations}\label{C}
\begin{align}
\widetilde{C}:=\const_{m,x_0,S,\widetilde{D},b}^+&:= \frac{b} {\inf\limits_{x\in \partial S}  g_{\widetilde{D}}(x, x_0)}>0,
\tag{\ref{C}c}\label{cz0C}
\\
\overline{C}_M:=	\int_{\widetilde{D}\setminus \{x_0\}} g_{\widetilde{D}}(\cdot, x_0)  \dd {\mu}_M  &+\int_{\widetilde{D}\setminus S} g_{\widetilde{D}}(\cdot, x_0)  \dd {\mu}_M^-  +M^+(x_0), 
\tag{\ref{C}M}\label{CM}
\end{align}
\end{subequations}
причем в случае   $\widetilde{D}\Subset D$ имеем
$\overline{C}_M\overset{\eqref{CM}}{=}\const_{m,x_0,S, \widetilde{D},M,D}^+<+\infty$.
\end{theoA}
Будем использовать следующую упрощенную версию теоремы A.
\begin{theoB} При соглашениях  \eqref{SD0D} и  \eqref{SD0D+} и с той же функцией 
$M{\in} \dsbh_*(D)$  с зарядом Рисса $\mu_M\in \Meas(D)$ для любой функции $u\in \sbh_*(D)$ с мерой Рисса  $\mu_u\in \Meas^+(D)$, удовлетворяющей неравенству $u \leq M$ на $D$, найдется постоянная  $C:=\const_{m, D, S,u,M}^+\in \RR^+$, с которой имеют место  неравенства   
\begin{equation}\label{mest}
\int_{D\setminus S}  v \dd {\mu}_u 		\leq	\int_{D\setminus S}  v \dd {\mu}_M	+C\quad
\text{для всех $v\overset{\eqref{<b}}{\in}  \sbh_0^+(D\setminus S;\leq 1) $}. 
\end{equation}
\end{theoB}
\begin{proof} Найдутся точка  $x_0\in \Int S$ и $r_0\in \RR^+_*$, для которых  \cite[3.1]{KhaRoz18} 
\begin{equation}\label{z0r0}
\begin{split}
D(x_0,r_0)\Subset \Int S, \quad u(x_0)&\neq -\infty, \quad M(x_0)\neq \pm\infty , 
\\ 
\Bigl|\int_{D(x_0,r_0)}h_m (|x-x_0|) \dd \mu_M\Bigr|<+\infty;\quad &h_m(t):=\begin{cases}
\log t &\text{ при $m=2$},\\
-|t|^{2-m}&\text{ при $m>2$}.
\end{cases}
\end{split}
\end{equation}
Найдется регулярная для задачи Дирихле область $\widetilde{D}$, для которой 
$\Int S\Subset \widetilde{D}\Subset D$ {\rm \cite[2.6.2-3]{HK}}. Выбор такой точки  $x_0$ 
и такой области $\widetilde{D} $ предопределен исключительно множествами $S,D$. Выберем  $b:=1$.
Таким образом, $\widetilde{C}\overset{\eqref{cz0C}}{=}\const_{m,D,S}^+\in \RR_*^+$ --- постоянная, зависящая только от $m,D,S$. Ввиду  \eqref{z0r0} по определению \eqref{CM} постоянная  
$\overline{C}_M\overset{\eqref{CM}, \eqref{z0r0}}{=}\const_{m, D,S,M}^+\in \RR^+$  зависит только от 
$m, D,S,M$. Отсюда постоянная 
\begin{equation*}
C\overset{\eqref{mest+}}{:=}|\widetilde{C} u(x_0)|	+ 
|\mu_M|(\widetilde{D}\setminus S)+\widetilde{C}\, \overline{C}_M\\
\geq - \widetilde{C} u(x_0) + 
\int_{\widetilde{D}\setminus S} v \,d {\mu}_M^-   +\widetilde{C}\, \overline{C}_M,
\end{equation*}
зависит только от $m, D,S, u,M$, т.\,е.  $C=\const_{D,S,u,M}^+\in \RR^+$. Таким образом, \eqref{mest} следует из \eqref{mest+}.
\end{proof}
Метод построения субгармонических тестовых функция на  $\BB$ вне  $r\BB$ на основе  $\rho$-субсферических положительных функций дает следующее 
\begin{propos}\label{trc-stf} Пусть  $h\overset{\eqref{SS+}}{\in}\rho\text{-}\sbs^+(\SS)$, $\SS=\SS_{m-1}\subset \RR^m$, и  $g\colon \RR^+\to \RR^+$ --- выпуклая функция с  $g(0)=0$. Положим 
\begin{equation}\label{rrho}
  \frac{1}{2}\leq r_{\rho}:=\max\left\{\frac{1}{2}, 
\sqrt[m-1]{\Bigl(1-\frac{m-1}{\rho(\rho+m-2)}\Bigr)^+}
\right\}<1.
\end{equation} 
Тогда функция
\begin{equation}\label{gh}
x\mapsto g\Bigl(\frac{1}{|x|^{m-1}}-1\Bigr) h\bigl(x/|x|\bigr), \quad x\in \BB\setminus \{0 \},
\end{equation}     
принадлежит классу\/  {\rm (см. \eqref{<b})}
\begin{equation}\label{hsbh}
\sbh_0^+\bigl(\BB \setminus  r_{\rho}\clos  \BB; \leq b_{\rho}\bigr), \text{ где }b_{\rho}:=g\Bigl(\frac{1}{r_\rho^{m-1}}-1\Bigr)\max_{\SS} h.
\end{equation}
\end{propos}
\begin{proof} Используем свойства \eqref{tr:ii}--\eqref{tr:6} функций класса
 $\rho$-$\sbs (\SS)$.

Существует убывающая последовательность выпуклых положительных функций  $g_n \underset{n\to\infty}{\searrow} g$ на  $\RR$, для которой  $g_n(0)=0$ и  $g_n\in C^{2}(\RR^+_*)$, $n\in \NN$. Существует также убывающая последовательность 
 $\rho$-субсферических положительных функций  $h_n\underset{n\to\infty}{\searrow} h$, $h_n\in C^2(\SS)$, $n\in \NN$ \cite[предложение  1.4]{Kha91}.
Предел каждой убывающей последовательности субгармонических положительных функций --- положительная субгармоническая функция.
Следовательно, достаточно установить субгармоничность функции \eqref{gh} на  $\BB$ вне  $r_{\rho}\BB$ в случае  $h\in C^2(\SS)$ и  $g\in C^2(\RR^+_*)$.  
Для оператора Лапласа $\Delta$  в обозначении из \eqref{Lrho} для сферической части $\Delta_{\SS}$ оператора Лапласа с  $r:= |x|$ 
и $s:=|x|/r$  в сферических координатах имеем 
\begin{equation}\label{DeltaS}
\Delta= \frac{\partial^2}{\partial r^2}+\frac{m-1}{r}\frac{\partial}{\partial r}
+\frac{1}{r^2} \Delta_{\SS}, 
\end{equation}
откуда вычисление оператора Лапласа от функции  \eqref{gh}  дает
 \begin{multline}\label{Delta}
\Delta \left(g\Bigl(\frac{1}{r^{m-1}}-1\Bigr)h(s)\right)\overset{\eqref{gh},\eqref{DeltaS}}{=}h(s)\Bigl(\frac{\partial^2}{\partial r^2}+\frac{m-1}{r}\frac{\partial}{\partial r}\Bigl) g\bigl(r^{1-m}-1\bigr)
\\+\frac{1}{r^2} g\bigl(r^{1-m}-1\bigr)\bigl(\Delta_{\SS} h\bigr)(s)
 \\
=\Bigl(\frac{(m-1)^2}{r^{2m}} g''\bigl(r^{1-m}-1\bigr)+
\frac{m-1}{r^{m+1}}g'\bigl(r^{1-m}-1\bigr)\Bigr)h(s)
+\frac{1}{r^2} g\bigl(r^{1-m}-1\bigr)\bigl(\Delta_{\SS} h\bigr)(s)
\\
\geq \frac{m-1}{r^{m+1}}g'\bigl(r^{1-m}-1\bigr)h(s)
+\frac{1}{r^2} g\bigl(r^{1-m}-1\bigr)\bigl(\Delta_{\SS} h\bigr)(s),
\end{multline}
так как по условию $h\geq 0$ на $\SS$, а также $g''\geq 0$ для выпуклой функции $g \colon \RR^+\to \RR^+$. Кроме того, при $g(0)=0$ каждая такая функция $g$ обладает свойствами 
\begin{equation}\label{gx}
 g'(x)\geq \frac{g(x)}{x} \quad\text{\it для всех\/ $x\in \RR_*^+$, \quad $g\in C(\RR^+)$ --- возрастающая.}
\end{equation}
Ввиду положительности функции $h$ из \eqref{Delta} и  \eqref{gx}  следует, что 
\begin{equation}\label{posgh}
\Delta \left(g\Bigl(\frac{1}{r^{m-1}}-1\Bigr)h(s)\right)
\geq \frac{1}{r^2}g\bigl(r^{1-m}-1\bigr)
\Bigl( \frac{m-1}{1-r^{m-1}}h(s)+\bigl(\Delta_{\SS} h\bigr)(s) \Bigr)
\end{equation}
{\it для  всех $r\in \RR_*^{+}$, $s \in \SS$.\/} По теореме S(s\ref{ss4}) для функции $h\in 
\rho\text{-}\sbs (\SS)$ из 
\eqref{Lrho} имеем 
\begin{equation}\label{DS}
\bigl(\Delta_{\SS} h\bigr)(s)\overset{\eqref{Lrho}}{=}(\mathcal L_{\rho} h)(s)- \rho (\rho+m-2) h(s)\overset{\eqref{Lrho}}{\geq} - \rho (\rho+m-2) h(s), \quad s\in \SS.
\end{equation}
Из последнего и из  \eqref{posgh} получаем 
\begin{equation}\label{posgh+}
\Delta \left(g\Bigl(\frac{1}{r^{m-1}}-1\Bigr)h(s)\right)
\overset{\eqref{DS}}{\geq} \frac{1}{r^2}g\bigl(r^{1-m}-1\bigr) h(s)
\Bigl( \frac{m-1}{1-r^{m-1}}-\rho(\rho+m-2)\Bigr).
\end{equation} 
Если $r\overset{\eqref{rrho}}{\geq} r_{\rho}$, то последняя скобка положительна, а следовательно  положительна  и  правая часть неравенства  
\eqref{posgh+}. Таким образом, функция  \eqref{gh} субгармоническая на  $\BB\setminus r_{\rho}\clos  \BB$. Очевидно, функция  \eqref{gh} положительна, так как   $h\in \rho\text{-}\sbs^+(\SS)$ и  $g\colon \RR^+\to \RR^+$ положительны. Кроме того, 
 \begin{equation}\label{g0}
 g (0) = 0 \quad \Longrightarrow \quad  \lim_{0<x\to 0}g(x)\overset{\eqref{gx}}{=}0 \quad \Longrightarrow \quad \lim_{1>r\to 1} g\Bigl(\frac{1-r}{r}\Bigr)h(\theta)\overset{\eqref{gh}}{=}0
\end{equation}
и, ввиду \eqref{gx}, имеем 
\begin{equation*}
g\Bigl(\frac{1}{r^{m-1}}-1\Bigr)\max_{\SS}h \overset{\eqref{rrho}}{\leq} 
g\Bigl(\frac{1}{r_{\rho}^{m-1}}-1\Bigr)\max_{\SS}h
\overset{\eqref{hsbh}}{=} b_{\rho} \quad\text{для всех  $r\in (r_{\rho}, 1)$}.
\end{equation*}
По определению \ref{def:tesf} согласно \eqref{g0} функция  \eqref{gh} принадлежит классу \eqref{hsbh}.  
\end{proof}

\subsection{Доказательства основных результатов} 

\begin{proof}[Основной теоремы] Пусть число $\rho\geq 0$ таково, что 
$r_{\rho}\overset{\eqref{rrho}}{=}1/2$.  
По предложению  \ref{trc-stf} функция   \eqref{gh} принадлежит классу 
  $\sbh_0^+\bigl(\BB \setminus \frac12\clos  \BB;\leq 1\bigr)$, так как  $g(1)\leq 1$ и  $\max_{\SS} h\leq 1$ по условиям  [g]--[h] основной теоремы. 
Отсюда по теореме B      существует постоянная  $C=\const_{m, M,u}^+$, для которой неравенство \eqref{mest} с $D=\BB$ и $S=\frac12\BB$  выполнено для любых функций  $v$ вида \eqref{gh}. 
Таким образом, получаем  
\begin{multline}\label{uM+}
\int_{\BB\setminus \frac12\clos \BB} g\Bigl(\frac{1}{|x|^{m-1}}-1\Bigr) h\bigl(x/|x|\bigr)
\dd \mu_u(x)
\\
\overset{\eqref{mest}}{\leq} \int_{\BB\setminus \frac12\clos \BB}  g\Bigl(\frac{1}{|x|^{m-1}}-1\Bigr) h\bigl(x/|x|\bigr) \dd \mu_M(x)+C
\quad\text{\it при\/  {\rm [g]}--{\rm [h]}}.
\end{multline}
\begin{lemma}[{\rm (\cite[\S~3]{Kha91}, \cite[(4.2)--(4.3)]{Kha99}, \cite[4.2.2]{Khsur})}]\label{l1} Пусть $r\in (0,1)$,  $f\in C(r,1)$,  т.\,е. функция $f$ непрерывна на интервале $(r,1)\subset \RR$,  $\mu \in \Meas \BB$. Тогда в обозначениях  \eqref{ntk} имеет место равенство 
\begin{equation}\label{uM+S}
\int_{\BB\setminus r\clos \BB} f\bigl(|x|\bigr) h\bigl(x/|x|\bigr)\dd \mu (x)\overset{\eqref{irR}}{=}
\int_{r}^1 f(t) \dd \mu^{\rad}(t;h).
\end{equation}
\end{lemma}
По лемме \ref{l1} из \eqref{uM+} и \eqref{uM+S}  следует  заключение  \eqref{uM} основной теоремы  при $r_{\rho}=1/2$.

Рассмотрим теперь случай  $1/2<r_{\rho}\overset{\eqref{rrho}}{<}1$.  По предложению  \ref{trc-stf} функция  \eqref{gh} принадлежит классу   \eqref{hsbh}. При этом  
$$\sbh_0^+\bigl(\BB \setminus r_{\rho}\clos  \BB;\leq b_{\rho}\bigr)\subset  \sbh_0^+\bigl(\BB \setminus r_{\rho}\clos  \BB;\leq 1\bigr),$$  
так как  $g(1)\leq 1$ и $\max_{\SS} h\leq 1$ при  ограничениях {[g]}--{[h]}, а также 
\begin{equation*}
b_{\rho}\overset{ \eqref{hsbh}}{\leq} g\Bigl(\frac{1}{(1/2)^{m-1}}-1\Bigr)\max_{\SS} h\overset{\eqref{gx}}{\leq} g(2^{m-1}-1)\max_{\SS} h\overset{{\rm [g]\text{--}[h]}}\leq 1 
\quad \text{ при $r_\rho>1/2$}.
\end{equation*}
Отсюда по теореме B существует постоянная  $C'=\const_{m, S,M,u}^+=\const_{m, \rho,M,u}^+$ для  $S:=r_{\rho}\clos \BB$, для которой неравенство  \eqref{mest} выполнено для любой функции  $v$ вида \eqref{gh},  а именно:
\begin{multline}\label{uM+r}
\int_{\BB\setminus r_\rho\clos \BB} g\Bigl(\frac{1}{|x|^{m-1}}-1\Bigr) h\bigl(x/|x|\bigr) \dd \mu_u(x)
\\ \overset{\eqref{mest}}{\leq} 
\int_{\BB\setminus r_\rho\clos \BB} g\Bigl(\frac{1}{|x|^{m-1}}-1\Bigr) h\bigl(x/|x|\bigr)\dd \mu_M(x)+C'
\quad\text{\it при\/ {\rm [g]}--{\rm [h]}}.
\end{multline}
Легко видеть, что существуют  $C'':=\const_{m, \rho, u}^+$, $C''':=\const_{m, \rho, M}^+$, для которых 
\begin{subequations}
\begin{align*}
\int_{r_\rho \clos \BB\setminus \frac12\clos \BB} g\Bigl(\frac{1}{|x|^{m-1}}-1\Bigr) h\bigl(x/|x|\bigr) \dd \mu_u(x)
\\
\leq \mu_u\bigl(r_\rho \clos \BB\setminus (1/2)\clos \BB\bigr)&\leq C'', \\
\left|\int_{r_\rho \clos \BB\setminus \frac12\clos \BB} g\Bigl(\frac{1}{|x|^{m-1}}-1\Bigr) h\bigl(x/|x|\bigr) \dd \mu_M(x)\right|
\\
\leq |\mu_M|\bigl(r_\rho \clos \BB\setminus (1/2)\clos \BB\bigr)&\leq C'''.
\end{align*}
\end{subequations}
Отсюда ввиду  \eqref{uM+r} получаем \eqref{uM+} с $C:=C'+C''+C'''= \const_{m, \rho,M,u}^+$
при выборе и ограничениях\/ {\rm [g]}--{\rm [h]}. 
По  равенству \eqref{uM+S} леммы \ref{l1} снова из  \eqref{uM+}  получаем заключение  \eqref{uM} основной теоремы уже для случая  $r_\rho>1/2$.

Пусть\/ ${\sf Z}$ --- поддивизор  нулей для функции $f\in \Hol_*(\BB)$,
$\BB\subset \CC^n$, удовлетворяющей неравенству $\log |f|\leq M$.\/ По заключению   \eqref{uM} основной теоремы  существует постоянная $C:=\const_{\rho,M,f}^+$, для которой имеем   \eqref{uM} с  $u:=\log |f|$ и с мерами Рисса $\mu_u=\mu_{\log|f|}$. Здесь выбор функции  $f$ предопределен исключительно функцией $\sf Z\colon \BB \to \RR^+$ и фукнцией   $M$, т.\,е. $C=\const_{\rho,M,{\sf Z}}^+$. 
Используя формулу Пуанкаре\,--\,Лелона, получаем цепочку  (не)равенств
\begin{multline*}
\int_{1/2}^{1} g\Bigl(\frac{1}{r^{2n-1}}-1\Bigr) {\dd} {\sf Z}^{\rad}(r ;h)
  \overset{\eqref{rcfZ}}{\leq} 
\int_{1/2}^{1} g\Bigl(\frac{1}{r^{2n-1}}-1\Bigr) \dd\, (\Zero_f)^{\rad}(r ;h)
\\
\overset{\eqref{nufZ}}{=}\int_{1/2}^{1} g\Bigl(\frac{1}{r^{2n-1}}-1\Bigr) \dd \mu_{\log |f|}^{\rad}(r ;h)
\leq \int_{1/2}^{1} g\Bigl(\frac{1}{r^{2n-1}}-1\Bigr) \dd \mu_{u}^{\rad}(r ;h)
\end{multline*}
для $u:=\log|f|$.  Таким образом, \eqref{uMf} следует из \eqref{uM}. 
\end{proof}

\begin{proof}[Теоремы единственности] Не умаляя общности, можем считать, что
 $h_0:=\max_{\SS}h>0$ и   $g(2^{2n-1}-1)>0$. 
Если  $f\in \Hol_*(\BB)$,  $|f|\leq \exp M$ на  $\BB$ и $\sf Z$ --- поддивизор нулей функции $f$, то по основной теореме с $m=2n$ имеем 
\begin{multline*}
\frac{1}{g(2^{2n-1}-1)h_0}\int_{1/2}^{1} g(1-r) {\dd} {\sf Z}^{\rad}(r ;h)
\\
\overset{\eqref{cuZ}}{\leq}
\frac{1}{g(2^{2n-1}-1)h_0}\int_{1/2}^{1} g\Bigl(\frac{1}{r^{2n-1}}-1\Bigr) {\dd} {\sf Z}^{\rad}(r ;h)
\\
= \int_{1/2}^{1} \frac{1}{g(2^{2n-1}-1)}\, g\Bigl(\frac{1}{r^{2n-1}}-1\Bigr) {\dd} {\sf Z}^{\rad}(r ;h/h_0)
\\
\overset{\eqref{uMf}}{\leq} 
\int_{1/2}^{1} \frac{1}{g(2^{2n-1}-1)}\, g\Bigl(\frac{1}{r^{2n-1}}-1\Bigr) {\dd} \mu_M^{\rad}(r ;h/h_0)+C
\\
\overset{\eqref{ntk}}{=}
\frac{1}{g(2^{2n-1}-1)h_0} \int_{1/2}^{1}  g\Bigl(\frac{1}{r^{2n-1}}-1\Bigr) {\dd} \mu_M^{\rad}(r ;h)+C
\\
\leq
\frac{1}{g(2^{2n-1}-1)h_0} \int_{1/2}^{1}  g\bigl(2^{2n-1}(2n-1)(1-r)\bigr) {\dd} \mu_M^{\rad}(r ;h)+C \overset{\eqref{1.n}}{<}+\infty.
\end{multline*}
Таким образом, если $f\neq 0$, то последнее противоречит условию  \eqref{cuZ}.
\end{proof}


\end{document}